\documentclass[12pt]{article}
\usepackage{amsmath}
\usepackage{amssymb}
\usepackage{amsthm}
\usepackage{amsfonts}
\usepackage{graphicx}
\usepackage{subfigure}
\usepackage{url}

\usepackage{bbm,epsfig,graphics,epic,color,rotating,color}
\numberwithin{equation}{section}

\newtheorem{theorem}{Theorem}[section]
\newtheorem{lemma}[theorem]{Lemma}
\newtheorem{proposition}[theorem]{Proposition}

\newtheorem{definition}[theorem]{Definition}

\theoremstyle{definition}

\title{Type-dependent stochastic Ising model describing the dynamics of a non-symmetric feedback module\thanks{Final version, as will appear in Mathematical Biosciences and Engineering}.}
\date{}
\author{Manuel Gonz\'alez-Navarrete\thanks{Instituto de Matem\'atica e Estat\'i{}stica, Universidade de S\~ao Paulo. Rua do Mat\~{a}o, 1010, 05508-090, S\~{a}o Paulo, Brazil. e-mail: manuelg@ime.usp.br}}

\begin{document}

\maketitle %
\thispagestyle{empty} %
\baselineskip=14pt

\vspace{8pt}

%
%

\begin{abstract}
We study an alternative approach to model the dynamical behaviors of biological feedback loop, that is, a type-dependent spin system, this class of stochastic models was introduced by Fern\'andez et. al \cite{FFJ}, and are useful since take account to inherent variability of gene expression.
We analyze a non-symmetric feedback module being an extension for the repressilator, the first synthetic biological oscillator, invented by Elowitz and Leibler \cite{E&L}. We consider a mean-field dynamics for a type-dependent Ising model, and then study the empirical-magnetization vector representing concentration of molecules. We apply a convergence result from stochastic jump processes to deterministic trajectories and present a bifurcation analysis for the associated dynamical system. We show that non-symmetric module under study can exhibit very rich behaviours, including the empirical oscillations described by repressilator.
\end{abstract}


\section{Introduction}
\label{sec:intro}

Elowitz and Leibler \cite{E&L} addressed the design and construction of a synthetic network providing the basic functionality of generating oscillations. Such functionality is essential to organize time-modulated biological functions like, for instance, the required periodic adjustment of an organism's physiology to the circadian rhythm \cite{E&N,Ran}. Their idea was to construct a negative feedback loop composed by three transcriptional repressors that are not part of any natural biological clock. The resulting oscillating module was called {\it repressilator} (see Figure~\ref{fig:rep} for illustration).

We remark, as done by Elowitz and Leibler \cite{E&L}, that the repressilator displayed noisy behavior, this fact must be related with stochastic fluctuations of its components. Indeed, Elowitz et al. \cite{El2} verified that gene expression is inherently variable, or noisy, due to random fluctuations in individual cells (see also, Shahrezaei and Swain \cite{SandS} and Swain et al. \cite{Swa}).

Progress in the understanding of this module, as well as others naturally occurring networks, and their associated control mechanisms demand the development of mathematical models that manage good balance between simplicity and usefulness. In the literature there exist several different approaches to this aim. For instance, in the case of the repressilator we could find works such as Chen and Aihara \cite{Ch1}, Chen et al. \cite{Ch2} and Wang et al. \cite{Wang1}.

This kind of interaction modules are widely known as feedback loops, which carry out specific functions in a cell, decomposing realistic cellular control networks \cite{Alon,Son}. The study of feedbacks loops have received extensive attention in recent years \cite{BFG,MG,Rad}. Particularly, the development of general frameworks can be found in Along \cite{Alon2}, Tyson and Novak \cite{T&N}, and especially in Wang et al. \cite{Wang2}. In general, we are interested in understanding for the dynamical behaviour of the concentration of molecules involved in the interactions.

In this work, we propose an alternative approach, which consists in the application of a class of interacting particle systems (see \cite{Lig}). The ideas are taken from the type-dependent stochastic spin system, proposed by Fern\'andez et al. \cite{FFJ}. Essentially, the main feature of this approach is taking account the inherent random variability in gene expression. The potential applications in general biological signaling networks are discussed in Fern\'andez et al. \cite{FFJ}.

As an illustration of our approach, we propose the study of a particular cyclic-interaction loop, that we call {\it non-symmetric clock module}, which is a simple extension of repressilator. Thus, to address this kind of qualitative approach, we propose the application of a mean-field type-dependent Ising model dynamics. Although the Ising model was initially studied to understand the physical phenomenon of ferromagnetism \cite{Ising}, nowadays this model represents a useful tool in different areas such as image processing, neural networks or earthquake dynamics \cite{D&Z,bio,geo1}, among others.

In our modeling scheme, the dynamics of the type-dependent Ising models is projected onto associated continuous-time jump processes, called \textit{density-profile processes}, which are random walks mapping the macroscopic evolution of the particle systems.

Fern\'andez et al. \cite{FFJ} showed that for arbitrary but fixed time intervals, in the limit of a very large number of particles (thermodynamics limit), the evolution of these jump processes converges to time dependent functions satisfying a correspondent deterministic dynamical system. We will include a simpler and straightforward proof of the convergence from the stochastic trajectories of the density-profile to deterministic paths ruled by non-linear differential equations. Our technique is based on the work of Ethier and Kurtz \cite{E&K}, who characterized a class of Markov jump processes called \textit{density dependent population processes} (see also Kurtz \cite{Ku2}).

Therefore, we use the convergence result to study the dynamical behaviours of our non-symmetric clock module. Particularly, we will analyse the influence of parameters in the behaviours of our model. The characterization of the evolution of the associated dynamical systems, that is a bifurcation analysis, completes our work.

We remark that the approach introduced by Fern\'andez et al. \cite{FFJ} allows us to analyse general classes of feedback loops. For instance, we can enumerate the class of feedback loops studied in \cite{L&W}, which are simple enough to be analysed in a theoretical framework, but admiting very rich dynamical behaviours. Moreover, all the loops have been found in transcriptional regulatory networks (Leite and Wang \cite{L&W}).

As we will see in section \ref{sec:dynamics} the approach proposes a non-reversible stochastic  microscopic dynamics based on interacting particle systems and statistical mechanics ideas. We used a mean-field version as a suitable approximation of other Ising-type interactions. In this sense, Mendon\c{c}a and de Oliveira \cite{MO} discussed properties of the stationary measure of a TDSIM with near neighbours interaction in the study of repressilator. Further investigations could focus in other feedback loops, exploiting applications and to better understanding of the type-dependent stochastic dynamics.

The paper is organized in the following way. In section \ref{sec:module} we introduce the non-symmetric clock module, which will be used as an illustration. Section \ref{sec:general} includes the definition of type-dependent spin models, and particularly the dynamics of the mean-field type-dependent stochastic Ising model. We construct the associated density-profile process and prove convergence to deterministic trajectories in section \ref{sec:perfil}. Section \ref{sec:bifur} includes a bifurcation analysis of the dynamical system associated to our non-symmetric clock module. Finally, section \ref{sec:final} contains some conclusions and further investigations.

\bigskip


\section{Feedback loops: the example of a non-symmetric clock module}
\label{sec:module}

In this section we introduce a specific feedback module, which will be used to explain the step by step of our modeling scheme. However, we remark that the stochastic approach in Section \ref{sec:general} can be applied for dynamical studies of general feedback loops.

\begin{figure}[htp]
\begin{center}
\subfigure[]{
\epsfig{file=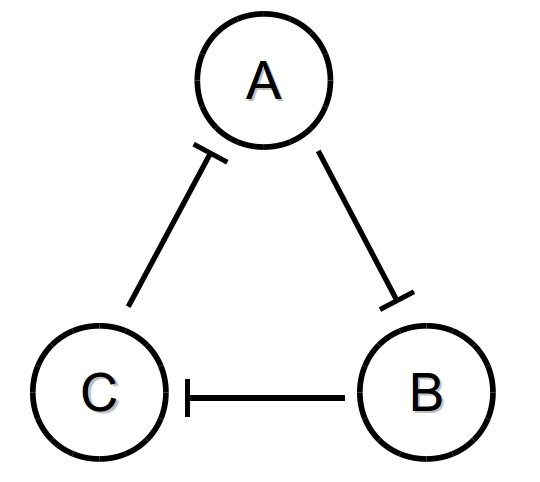, height=4cm}
\label{fig:rep}
}
\hspace{0.6cm}
\subfigure[]{
\epsfig{file=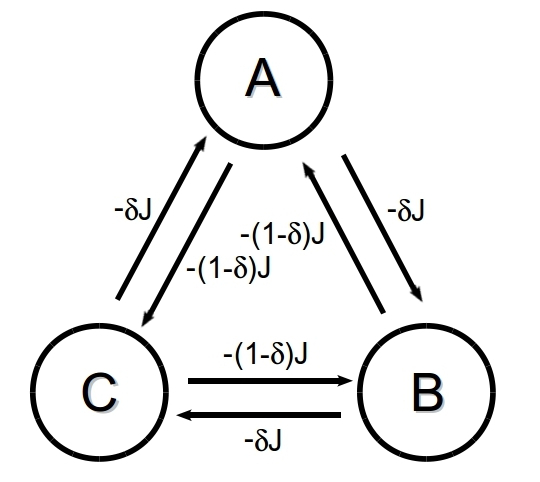, height=4cm}
\label{fig:nonsym}
}
\end{center}
\caption{{\scriptsize Representation of cyclic feedback modules: (a) Repressilator: a simple model of transcriptional regulation representing the loop of feedback inhibition among three components indicated by \textit{A}, \textit{B} and \textit{C}. The blunt arrows indicate inhibition. (b) Motif of our non-symmetric clock module (including the rates of interaction). The parameter $\delta \in [0, 1]$, and for positive or negative values of $J$, we have inhibition or activation cycles, respectively.}}\label{falgo}
\end{figure}

We use the notion of network motifs to design the repressilator and our non-symmetric clock module. A feedback loop motif is a simple representation of a transcriptional regulation network \cite{L&W,T&N}. That is a cycle in a directed graph whose vertices (also called nodes) can represent concentrations of proteins or genes (see Figure \ref{falgo}). In this sense, the edges can represent either positive or negative interactions. In other words, a positive (resp. negative) edge implies that the component in tail vertex activates (resp. represses) the transcription of the gene of the element in head vertex.

We shall consider a simple feedback loop that we call {\it non-symmetric clock module}, it is based on the repressilator proposed by Elowitz and Leibler \cite{E&L}. The repressilator is a three transcriptional repressor systems that are not part of any natural biological clock, and were used to build an oscillating network in Escherichia coli. In other words, this is a negative feedback loop which provides the basic functionality of generating oscillations. Which is biologically required in several contexts like in cell-cycle and in the setting up of circadian cycle.

Figure \ref{fig:rep} shows the representation of the repressilator. We denote by \textit{A, B} and \textit{C}, the \textit{LacI} protein, the \textit{TetR} gene and the \textit{CI} gene, respectively. Accordingly, the first repressor protein, LacI from E. coli, inhibits the transcription of the second repressor gene, tetR from the tetracycline-resistance transposon Tn10, whose protein product in turn inhibits the expression of a third gene, cI from l phage. Finally, CI inhibits lacI expression, completing the cycle. In other words, the rate of change of component \textit{A} density at each time depends only on \textit{C} density in an inhibitory manner: the density of \textit{A} tends to decrease if the concentration of \textit{C} is high. A similar dependence holds between \textit{C} and \textit{B} and between \textit{B} and \textit{A}.

In this work, we propose the study of a non-symmetric cyclic-interaction module, see Figure \ref{fig:nonsym}. We will refer to $A,B$ and $C$ as being three types of {\it molecules}. We focus on this simple feedback loop because it is large enough to admit complex dynamical behaviours. The model includes the parameter $J$ that measures the strength of the interaction and a parameter $\delta \in [0, 1]$, which allows us to distribute the interactions between each pair of the three components ($A,B,C$). Note that the interactions between neighbour components could be non-symmetric, but the global interaction holds the invariance.

The particular cases when $\delta$ is equal $0$ or $1$, are similar to the repressilator above defined. Although, it is important to stress that for positive or negative values of $J$, we have inhibition or activation cycles, respectively. We are interested in a characterization of dynamical behaviour for this module, as a function of parameter $J$ and $\delta$.



Therefore, we describe a microscopic stochastic approach to derive an associated dynamical systems, which only seeks to incorporate the essential qualitative information about the biochemical interactions. This approach is based in the work of Fern\'andez et al. \cite{FFJ}, which borrows ideas from interacting particle systems \cite{Lig}, in such a way that spins represent the internal states of components of the feedback loop.

\bigskip

\section{Modeling setup: microscopic type-dependent dynamics}
\label{sec:general}

In this section we describe the type-dependent stochastic spin systems proposed by Fern\'andez et al. \cite{FFJ} to study signaling biological networks. We focus our explanation by studying the dynamical behaviours of the non-symmetric clock module exposed in previous section. In particular, we will define a {\it type-dependent stochastic Ising model} (TDSIM). That is, a microscopic model with Ising-type interactions and having dynamical evolution defined by the type-dependent dynamics proposed in Fern\'andez et al. \cite{FFJ}.

\bigskip

\subsection{A family of interacting particle systems.}

The \textit{type-dependent stochastic spin models} are a family of stochastic spin-flip systems introduced in Fern\'andez et al. \cite{FFJ}. This family extends the usual definition of particle systems \cite{Lig}, to allow asymmetric dependence of rates on the energy function, the Hamiltonian. As a consequence of this asymmetric dependence, a particularity of these models is its non-reversible stochastic dynamics.

\begin{figure}[ht]
\centering
\includegraphics[scale=0.15]{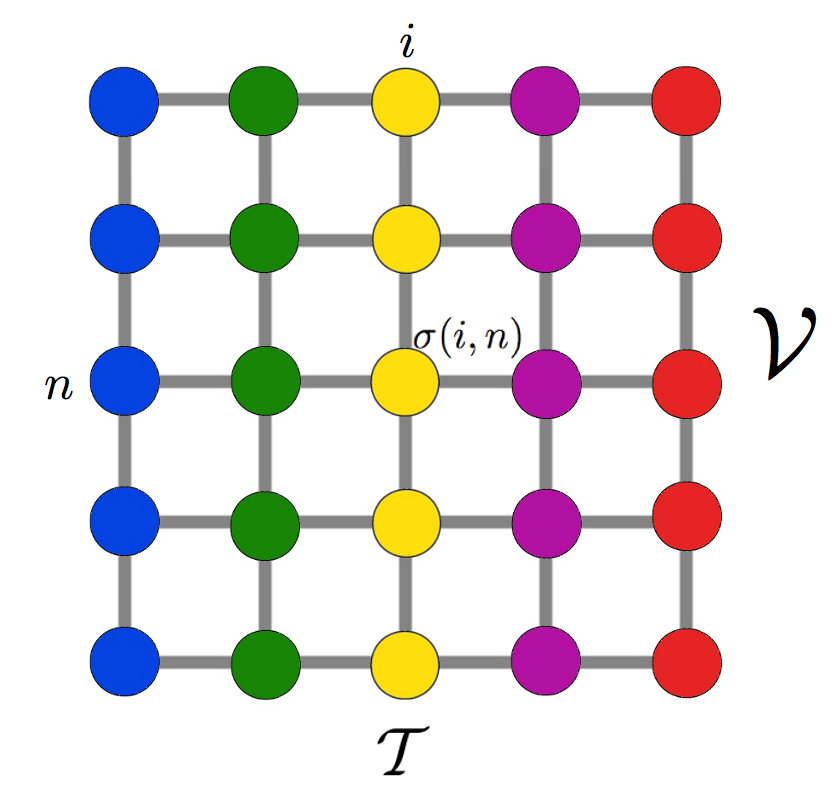}
\caption{{\scriptsize Lattice representation of the site-space $\Lambda$ for a type-dependent spin system. Each vertical line represents a type $i \in \mathcal{T}$. At horizontal lines we see the spatial positions $\mathcal{V}$. The central spin shows the notation $\sigma(i,n) \in \mathcal{S}_i$, for a given configuration $\sigma \in \Omega_{\Lambda}$.}}
\label{fig:model}
\end{figure}

Next, we explain these models. We need to introduce a set of {\it spin types} $\mathcal{T}$, of cardinal $k$, representing genes or proteins (for instance: LacI, TetR and CI, in repressilator above). Also, a vertex set $\mathcal{V}$ of a simple finite graph of order $N = | \mathcal{V}|$, denoting {\it spatial positions} available for each one of the types. We call the ordered pair $(i,n) \in \Lambda = \mathcal{T} \times \mathcal{V}$ a site. Moreover, the set of {\it internal states} for each $i \in \mathcal{T}$, will be denoted by $\mathcal{S}_i = \{a_1, \ldots, a_{s_i} \}$. Hence, the spin system has site-space $\Lambda$ and configuration space $\Omega_{\Lambda} = \prod_{i \in \mathcal{T}} \mathcal{S}_i^{\mathcal{V}}$. For a configuration $\sigma \in \Omega_{\Lambda}$ we denote by $\sigma(i,n)$ the value of the spin at site $(i, n)$ (see Figure \ref{fig:model}).\\

\bigskip 

The continuous-time evolution of these models is governed by a non-reversible Glauber spin-flip stochastic dynamics.  In other words, we have a stochastic evolution in continuous time, for which only one particle flips at each transition. Moreover, the dynamics is non-reversible with respect to the Gibbs measure (defined in \eqref{Gibbs}).

The corresponding rates are determined in terms of a function $H_{\Lambda}: \Omega_{\Lambda} \to \mathbb{R}$, the Hamiltonian of the spin system. To define the Hamiltonian, we denote
 
 \begin{equation}
\label{Cepsilon}
\mathcal{C} = \left\{ \underline{i} = (i,a): i \in \mathcal{T}, a \in \mathcal{S}_i \right\},
\end{equation}
\vspace{0.1cm}
 
\noindent then, the Hamiltonian of these models is determined by a family of interaction matrices $\mathbb{J}_{n,l} : \mathcal{C} \times \mathcal{C} \to \mathbb{R}$, one for each pair of spatial positions $n, l \in \mathcal{V} $. Due to the original applications in signaling biological networks, these matrices $\mathbb{J}_{n,l}[\cdot ; \cdot]$ are not assumed to be symmetric. Particularly, we say that $\mathbb{J}_{n,l} [(i,a);(j,b)]$ indicates the strength of the influence that a spin at a site $(i, n) \in \Lambda$ in internal state $a \in \mathcal{S}_i$ has upon a spin at $(j, l) \in \Lambda$ that is in internal state $b \in \mathcal{S}_j$. 

As we illustrated in the repressilator module (see Figure \ref{fig:rep}), type A components act upon type B components, while the reciprocal interaction does not occur. Then, it is reasonable that the quantity $\mathbb{J}_{n,l} [(i,a);(j,b)]$ may be noticeably different that the one in the opposite direction $\mathbb{J}_{l,n} [(j,b);(i,a)]$. However, we are not assuming asymmetry with respect spatial positions, because the interaction is only associated to types, that is

\begin{equation}
\label{assimetria}
\mathbb{J}_{n,l} \left[ (i, a); (j,b) \right] = \mathbb{J}_{l,n} \left[ (i, a); (j,b) \right],
\end{equation}
\vspace{0.1cm}

\noindent for $n,l \in  \mathcal{V}$, $i,j \in \mathcal{T}$, $a \in \mathcal{S}_i$ and $b \in \mathcal{S}_j$. Consequently, as usual in statistical mechanics, the Hamiltonian is then defined by

\begin{equation}
\label{hamiltoniano}
H_{\Lambda} (\sigma)= - \sum_{(i,n) \in \Lambda} \sum_{(j,l) \in \Lambda} \mathbb{J}_{l,n} \left[ (j, \sigma(j,l)); (i,\sigma(i,n)) \right],
\end{equation}
for each configuration $\sigma \in \Omega_{\Lambda}$.

\bigskip

\subsection{Type-dependent stochastic dynamics.}
\label{sec:dynamics}
We define the stochastic dynamics proposed by Fern\'andez et al. \cite{FFJ}. This dynamics allows only single spin-flips, more precisely, single-site internal-state transitions. In other words, we assume that at each time only one molecule could be produced or degraded. In statistical mechanics notation, we say that we have a Glauber type stochastic dynamics.

Formally, given a configuration $\sigma \in \Omega_{\Lambda}$, a site $(i,n)$ and an internal state $a \in \mathcal{S}_i$, we denote $\sigma^a_{(i,n)}$ the configuration with

\begin{equation}
\label{conf.estado}
[\sigma^a_{(i,n)}](j,l) = 
\left\{
\begin{array}{rl}  a, & \text{ if } (j,l)=(i,n), \\
\sigma(j,l), & \text{ otherwise}.
\end{array}
\right.
\end{equation}
\vspace{0.1cm}

That is, a configuration for which we fix the value $a$ for the spin at site $(i,n)$. Thus, the energy cost for the transition from $\sigma^a_{(i,n)}$ to $\sigma^b_{(i,n)}$ is given by

\begin{equation}
\label{energycost}
\Delta_{(i,n)}^{a \to b} (\sigma) =H(\sigma^b_{(i,n)}) - H(\sigma^a_{(i,n)}).
\end{equation}
\vspace{0.1cm}

Usually, this total change of the energy associated to the flip at site $(i,n)$ from state $a$ to $b$, is used to define the rates of transition in particle systems. Hence, generally we have a continuous-time stochastic evolution being reversible with respect to the Gibbs measure,

\begin{equation}
\label{Gibbs}
\mu_{\Lambda}(\sigma) = \displaystyle\frac{e^{-H_{\Lambda}(\sigma)}}{\sum_{\eta\in\Omega_{\Lambda}} e^{-H_{\Lambda}(\eta)}},
\end{equation}
for each $\sigma \in \Omega_{\Lambda}$.

However, in type-dependent spin models the definition of the rates of transition is quite different. The asymmetry of the interaction defined above leads naturally to the decomposition of \eqref{energycost} in the following manner,

 \begin{equation}
\label{energycost2}
\Delta_{(i,n)}^{a \to b} (\sigma) =\Delta[\text{IN}]_{(i,n)}^{a \to b} (\sigma) + \Delta[ \text{OUT}]_{(i,n)}^{a \to b} (\sigma),
\end{equation}
\vspace{0.1cm}

\noindent where, 
 
\begin{equation}
\label{sobre.ele}
\Delta[ \text{IN}]_{(i,n)}^{a \to b} (\sigma) = \displaystyle\sum_{(j,l) \in \Lambda} \left( \mathbb{J}_{l,n} \left[ (j,\sigma(j,l)); (i,a) \right] - \mathbb{J}_{l,n} \left[ (j,\sigma(j,l)); (i,b) \right] \right),
\end{equation}
\vspace{0.1cm}

\noindent that collects the change in the influence of the configuration $\sigma$ upon the site $(i,n)$, when internal state there changes from $a$ to $b$. On the other hand,

\begin{equation}
\label{sobre.outros}
\Delta[\text{OUT}]_{(i,n)}^{a \to b} (\sigma) = \displaystyle\sum_{(j,l) \in \Lambda} \left( \mathbb{J}_{n,l} \left[ (i,a) ; (j,\sigma(j,l)) \right] - \mathbb{J}_{n,l} \left[ (i,b); (j,\sigma(j,l)) \right] \right),
\end{equation}
\vspace{0.1cm}

\noindent collects the change of the influence that the site $(i,n)$ has on all other sites when its internal state flips from $a$ to $b$.

Furthermore, as a particularity of the type-dependent spin models, each transition rate depends only on the energy changes brought upon the site \eqref{sobre.ele}. Thus, we denote $\lambda^{a \to b}_{(i,n)}(\sigma)$ the rate of a transition flipping $\sigma^a_{(i,n)}$ to $\sigma^b_{(i,n)}$, which depends only on \eqref{sobre.ele}, that is,

\begin{equation}
\label{taxa.spin}
\lambda^{a \to b}_{(i,n)}(\sigma) = \Phi \left(\Delta[\text{IN}]_{(i,n)}^{a \to b} (\sigma)\right),
\end{equation}
\vspace{0.1cm}

\noindent where $\Phi$ is a non-increasing $\mathbb{R}_+$-valued function satisfying the \textit{detailed-balance} condition $\Phi(\Delta) e^{\Delta} =\Phi(-\Delta) e^{-\Delta}$.

Finally, we state a formal definition of the spin models proposed by Fern\'andez et al. \cite{FFJ}.

\medskip

\begin{definition}
\label{typedependent}
A {\it type-dependent stochastic spin model} is the continuous-time process $(\sigma_t)_{t \ge 0}$ on $\Omega_{\Lambda}$, defined by a spin model with a type-dependent interaction, given by the Hamiltonian in \eqref{hamiltoniano} and a dynamics with rates of the form \eqref{taxa.spin}.
\end{definition}

\medskip

The dynamics of these spin systems yields non-reversibility with respect to the Gibbs measure \eqref{Gibbs}. Then, the type-dependent dynamics has mathematical interest for itself. It is natural to study the stationary measure for particular cases of the interactions matrices $( \mathbb{J}_{\underline{i},\underline{j}})_{\underline{i},\underline{j} \in \mathcal{C}}$. That is, for instance local interactions like nearest neighbors (n.n), as usual in many particles systems. In Mendon\c{c}a and de Oliveira \cite{MO}, was proposed an Ising-type n.n interaction to analyse, by Monte Carlo simulations, some properties of the stationary measure for the repressilator.

In this work, as done in Fern\'andez et al. \cite{FFJ}, we consider a mean-field interaction for an Ising model. The mean-field interaction is a natural approximation of local interactions and given our interest in concentrations of molecule, as global observable, this is a suitable choice. Of course, the convergence results to be explained in Section \ref{sec:perfil} hold whenever we assume mean-field dynamics.

\medskip

\begin{definition}
\label{mean-field}
A type-dependent stochastic spin model is {\it mean-field} if the Hamiltonian parameters are of the form

\begin{equation}
\label{meanfield}
\mathbb{J}_{n,l}[(i,a);(j,b)] = \frac{\alpha_{(i,a),(j,b)}}{|\mathcal{V}|},
\end{equation}

\noindent where $\{ \alpha_{\underline{i},\underline{j}}\}_{\underline{i},\underline{j} \in \mathcal{C}}$, is a real matrix.
\end{definition}
\medskip

We conclude this section explaining the particular modeling setup that we will use to study the non-symmetric clock module.

\bigskip

\subsection{Mean-field TDSIM}
\label{sec:ising}

We follow the ideas of Fern\'andez et al. \cite{FFJ} to model our non-symmetric clock module through a \textit{type-dependent stochastic Ising model (TDSIM)}, which has all internal state spaces $\mathcal{S}_i=\{-1,+1\}$, for all $i \in \mathcal{T}$.

We think the spin types as points $\{A, B, C\}$ over the circle (see Figure \ref{fig:nonsym}) and, for each $i \in \mathcal{T} = \{A, B, C\}$, we denote $h(i)$ the neighbour of $i$ in the clockwise direction. And $a(i)$ to be the neighbour in the anti-clockwise direction. Borrowing statistical mechanical nomenclature, we say that for $J < 0 $, the activation interactions in the Figure \ref{fig:nonsym} are ferromagnetic. Moreover, for $J > 0$, the inhibition cycle is said to be an anti-ferromagnetic model.

Next, observe that the set \eqref{Cepsilon} will be particularly defined by

 \begin{equation}
\label{CIsing}
\mathcal{C}^* = \left\{ (A,+1); (A,-1);(B,+1);(B,-1);(C,+1);(C,-1) \right\}.
\end{equation}

In addition, the set of vertices $\mathcal{V}$ will represent a reservoir of capacity $N$ (actually, we have three reservoirs, one for each type of molecule in $\mathcal{T}$). In this way, for a configuration $\sigma \in  \{-1,+1\}^{\Lambda}$, the spin value $\sigma(i,n)=+1$ means that there is a molecule of type $i \in \mathcal{T}$ at spatial position $n \in \mathcal{V}$. Otherwise, if $\sigma(i,n)=-1$, it means the absence of molecules of type $i \in \mathcal{T}$ in the corresponding reservoir at spatial position $n \in \mathcal{V}$.

Acordingly, the mean-field TDSIM must be defined by rates functions as in \eqref{taxa.spin}, where $\Phi (\Delta ) = e^{- \Delta}$ (for details, see \cite{FFJ}). In addition, the real matrix $\{ \alpha_{\underline{i},\underline{j}}\}_{\underline{i},\underline{j} \in \mathcal{C}^*}$ as in \eqref{meanfield} is given by

\begin{equation}
\label{meu.alfa}
\alpha_{(i,a),(j,b)} = \left\{\begin{array}{rl}  - \delta J b, & \text{ if } j = h(i)\text{ and } a=+1,\\  - (1- \delta) J b, & \text{ if } j = a(i)\text{ and } a=+1,\\ \kappa_i, & \text{ if } j = i,  \\    0,  & \text{ otherwise, }\end{array} \right.
\end{equation}
\vspace{0.1cm}
\noindent for $(i,a), (j,b) \in \mathcal{C}^*$, and where $J \in \mathbb{R}$, $\delta \in [0,1]$, and $\kappa_i \in \mathbb{R}$ is an external field (or chemical potential). It must be clear that a site $(i,n) \in \Lambda$ in internal state $-1$, will not influence others sites, because at $(i,n)$ there is no molecule.

Notice that, as we said above, for $\delta = 0$ (as well as $\delta=1$) and $J > 0$, the totally asymmetric case, we obtain an Ising model for the repressilator represented in Figure \ref{fig:rep}.

Finally, the transition rates for TDSIM are defined for each $\sigma \in \Omega_{\Lambda}$ by
\begin{equation}
\label{taxas.ising}
\begin{array}{l}
\lambda^{-1 \to +1}_{(i,n)}(\sigma) = \exp \left( + 2\left[-\frac{ \delta J}{|\mathcal{V}|}\cdot a^+_i -\frac{(1- \delta) J}{|\mathcal{V}|}\cdot h^+_i  + \kappa_i \right]\right), \ \text{and}\\[0.5cm]
\lambda^{+1 \to -1}_{(i,n)}(\sigma) =  \exp \left( - 2\left[-\frac{ \delta J}{|\mathcal{V}|}\cdot a^+_i -\frac{(1- \delta) J}{|\mathcal{V}|}\cdot h^+_i  + \kappa_i \right]\right),
\end{array}
\end{equation}
\vspace{0.1cm}
where  $a^+_i = |  \{ l \in \mathcal{V} : \sigma(a(i),l) = +1  \}|$, and $h^+_i=|  \{ l \in \mathcal{V} : \sigma(h(i),l) = +1  \}|$, for $i \in \{A,B,C\}$.  For each type $i$, the parameter $\kappa_i$ must be interpreted as the own rate of production of molecules $i \in \mathcal{T}$. In general $\kappa_i$ is a linear function of parameter $J$. Because we are interested in the relation between two-body interaction ($J$) with individual interaction ($\kappa$).

We summarize our modeling setup in the following statement:

\medskip

\begin{definition}
\label{meanfieldTDSIM}
The microscopic evolution of the non-symmetric clock module is described by a mean-field type-dependent stochastic Ising model with continuous-time Glauber dynamics $(\sigma_t)_{t\ge0}$, whose rates of transition are given by \eqref{taxas.ising}.
\end{definition}

\medskip

In the next section we will focus in the macroscopic evolution of the non-symmetric clock module. Thus, we shall define the density-profile processes, which are a family of random walks with continuous time evolution. Moreover, we will prove an almost sure convergence from these stochastic trajectories to deterministic paths governed by non-linear differential equations.

\bigskip


 \section{Macroscopic evolution: the density-profile process and its convergence to deterministic dynamics}
\label{sec:perfil}

In this section we focus our attention in the concentrations of biochemical components into our feedback loop defined in Section \ref{sec:module}. Hence, we will mostly be interested on the vector of empirical magnetization of the mean-field TDSIM. Finally, we will study the thermodynamic limit, that is, the behaviours when $|\mathcal{V}| = N \to \infty$.

We project the mean-field TDSIM onto a jump Markov process in $\mathbb{R}^3$, called density-profile process, that represents the densities of the three biochemical components in our feedback loop. Thus, in statistical mechanics notation, the stochastic spin system and the density-profile process provide, respectively, the microscopic and macroscopic views of the same model.

\medskip
\begin{definition}
\label{def1}
The \textit{density-profile process} associated to the non-symmetric clock module is a continuous-time jump process $(X^N(t))_{t\ge0} \in \mathcal{D}_N = \left\{[0,1] \cap \{k/N, k \in \mathbb{Z}\}\right\}^3$, for $|\mathcal{V}|=N \ge 1$. That is defined for each $t\ge0$ by

\begin{equation}
\label{densityp}
X^N(t) : = (X_A^N(t),X_B^N(t),X_C^N(t)), 
\end{equation}
where

\begin{equation}
\label{densityp2}
X^N_i(t)  = \displaystyle\frac{\left|  \{ n \in \mathcal{V} : \sigma_t(i,n) = +1  \} \right|}{N},
\end{equation}
\vspace{0.1cm}
\noindent $i \in \{A,B,C\}$. Moreover, the jumps of this process are of the form $\mathbf{l}/N$, where $\mathbf{l}$ belong to $\mathcal{J}=\{\pm(1,0,0); \pm(0,1,0);\pm(0,0,1)\}$.
\end{definition}

\medskip

We describe the evolution of that processes in the following way: for each position $X^N(t) = x \in \mathcal{D}_N$, the jumps of the form

\begin{equation}
\label{salto}
x \rightarrow x + \frac{\mathbf{l}}{N},
\end{equation}
where $\mathbf{l} \in \mathcal{J}$, are defined by a collection $\beta_{\mathbf{l}}(x)$ of functions, $\beta_{\mathbf{l}}: \mathcal{D}_N \to [0,\infty], \mathbf{l} \in \mathcal{J}$. Of course, we require that $x \in \mathcal{D}_N$ and $\beta_{\mathbf{l}}(x)>0$, imply $x+\mathbf{l}/N \in \mathcal{D}_N$. It must be clear that each $\beta_{\mathbf{l}}(x)$ will be defined as a function of the transition rates in \eqref{taxas.ising}.

Therefore, the jump in \eqref{salto} occurs with rate

\begin{equation}
\label{taxa}
N\beta_{\mathbf{l}}(x) = \frac{d}{ds} P\left( X^N(t+s) = x + \frac{\mathbf{l}}{N} \right) \bigg|_{s=0},
\end{equation}
where

\begin{equation}
\label{beta}
\beta_{\mathbf{l}}(x) =
\left\{
\begin{array}{rl}  x_i \lambda^{+1 \to -1}_{(i,n)}(\sigma_t), & \text{ if  } \mathbf{l} = -e_i, \\ [0,3cm]
(1-x_i) \lambda^{-1 \to +1}_{(i,n)}(\sigma_t), & \text{ if  } \mathbf{l} = e_i,
\end{array}
\right.
\end{equation}
where $x_i=X_i^N(t)$ as defined in \eqref{densityp2}, and $e_i$ is the unitary vector in the direction of $i \in \{A,B,C\}$.

\vspace{0.1cm}

Note that each variable $x_i$ represents the density of elements of type $i$, or we say as well, the concentration of the component $i\in \{A,B,C\}$ in our non-symmetric clock module. Thus, $Nx_i$ will represent the number of molecules of the type $i$, and in some abstract sense, $N(1-x_i)$ denotes the available number of possible molecules into the reservoir of size $N$. Although, in the Ising system, this last quantity is just the number of sites of type $i \in \mathcal{T}$, with spin-value equal to $-1$.

The next step is the characterization of these jump processes in the thermodynamics limit, that is $N \to \infty$. Particularly, Fern\'andez et al. \cite{FFJ}, showed that in the limit of a very large number of spins, these density-profile processes converge to time dependent functions satisfying a correspondent deterministic dynamical system. They used a pathwise approach, which strongly exploits large deviation theory \cite{OV}, and coupling of random variables \cite{Coupling}. Their method provides explicit bounds for the distance between the stochastic and deterministic trajectories.

However, the qualitative study of the behaviors of our feedback loop, does not require this kind of accurate results obtained in \cite{FFJ}. Therefore, we follow the works of Kurtz \cite{Ku2} and Ethier and Kurtz \cite{E&K}, to obtain a simpler and straightforward proof of the convergence from the paths of the density-profile processes to deterministic trajectories, these latter ruled by non-linear differential equations.

Now, we aim to state our first theorem and prove it. The main tool that we will use is the characterization of a family of jump processes given by Kurtz \cite{Ku2}. We stress that our result can be easily extended to more general models, that include applications to a very extensive class of feedback loops.

Initially, we define a jump Markov process $\widehat{X}^N$ on $\mathbb{Z}^3$, for which we will use some properties given by Kurtz \cite{Ku2}. Particularly, the jumps of this process are of the form $\mathbf{l} \in \mathcal{J}$ (see Definition \ref{def1}), and their intensities will be given by $ \beta^*_{\mathbf{l}} (k) = N \beta_{\mathbf{l}} (k/N)$, $k \in \mathbb{Z}^3$, where $\beta_{\mathbf{l}}(x)$ as defined in \eqref{beta}. Thus, let $ \widehat{X}^N (0)$ being non-random, the evolution of $\widehat{X}^N$ at time $t\ge0$, can be written as

\begin{equation}
\label{process1}
\widehat{X}^N (t) = \widehat{X}^N ( 0) + \sum_{\mathbf{l}} \mathbf{l} Z_{\mathbf{l}}(t),
\end{equation}
\vspace{0.1cm}

\noindent where $Z_{\mathbf{l}}(t) : = | \{s \le t; \widehat{X}^N ( s) - \widehat{X}^N ( s-) = \mathbf{l}\}|$, that is, the number of jumps of size $\mathbf{l}$ until time $t$. Then, by some results from Chapter 7 of \cite{Ku2}, we have:

\begin{equation}
\label{process2}
\widehat{X}^N ( t) = \widehat{X}^N ( 0) + \sum_{\mathbf{l}} \mathbf{l} Y_{\mathbf{l}}\left(N\displaystyle \int_0^t \beta_{\mathbf{l}} \left(\frac{\widehat{X}^N (s)}{N}\right)ds\right),
\end{equation}
\vspace{0.1cm}

\noindent for each $t\ge0$, where the $Y_{\mathbf{l}}(u)$, $\mathbf{l} \in \mathcal{J}$, represent independent Poisson processes with corresponding intensities $u$.

Now, setting

\begin{equation}
\label{F(x)}
F(x) = \sum_{\mathbf{l}} \mathbf{l} \beta_{\mathbf{l}}(x),
\end{equation}
and $X^N = N^{-1}\widehat{X}^N$. Thus, we have

\begin{equation}
\label{process3}
X^N (t) = X^N ( 0) + \sum_{\mathbf{l}} \frac{\mathbf{l}}{N} \widetilde{Y}_{\mathbf{l}}\left(N\displaystyle \int_0^t \beta_{\mathbf{l}} \left(X^N(s)\right)ds\right) + \displaystyle \int_0^t F(X^N(s))ds,
\end{equation}
\vspace{0.1cm}

\noindent for each $t\ge0$, where $\widetilde{Y}_{\mathbf{l}}(u) = Y_{\mathbf{l}}(u) - u$, is the Poisson process centred at its expectation. Thus, our first result is the almost sure convergence of stochastic trajectories \eqref{densityp} to associated deterministic paths.

\bigskip


\begin{theorem}
\label{th1}

Consider the density-profile process $X^N(t)$, as defined in \eqref{densityp}-\eqref{densityp2}, with intensities given by \eqref{taxa}, \eqref{beta} and \eqref{taxas.ising}, which satisfies \eqref{process3}. Suppose $\lim_{N \to \infty} X^N ( 0 ) = x_0$, and $X$ satisfying

\begin{equation}
\label{condition3}
X(t) = x_0 + \displaystyle\int_0^t F(X(s))ds, \ \ t \ge 0.
\end{equation}
\vspace{0.1cm}

Then for every $t \ge 0$,

\begin{equation}
\label{result1}
\displaystyle\lim_{N \to \infty} \sup_{s \le t} |X^N(s)- X(s)| = 0 \ \ a.s.
\end{equation}

%
\end{theorem}
\medskip

The proof of this result is essentially repeat the arguments of Theorem 2.1 from Chapter 11 of \cite{E&K}, but we include the proof for completeness.

\begin{proof}

First of all, denote $\bar{\beta}_{\mathbf{l}} := \sup_{x \in \mathcal{D}_N} \beta_{\mathbf{l}} (x)$. Thus, it is easy to see that

\begin{equation}
\label{betafinite}
\sum_{\mathbf{l}} | \mathbf{l} | \bar{\beta}_{\mathbf{l}} < \infty,
\end{equation}

\noindent and by the Lipschitz continuity of the rates of the jumps \eqref{beta}, there exists a fixed $M > 0 $ such that

\begin{equation}
\label{F2}
| F(x) - F(y) | \le M | x- y |, \ \ x,y \in \mathcal{D}_N.
\end{equation}
\vspace{0.1cm}

Observe that by \eqref{process3} and \eqref{condition3}, we have for each $t\ge0$,

\begin{equation}
\label{distance1}
\begin{array}{lll}
|X^N(t) - X(t) | &\le& |X^N(0) - x_0 | +\left| \displaystyle\sum_{\mathbf{l}} \frac{\mathbf{l}}{N} \widetilde{Y}_{\mathbf{l}}\left(N\displaystyle \int_0^t \beta_{\mathbf{l}} \left(X^N(s)\right)ds\right) \right| \\[0.4cm]
&&+ \left| \displaystyle\int_0^t [F(X^N(s)) -F( X(s)) ]ds \right|.
\end{array}
\end{equation}

Now, let denote

\begin{equation}
\label{epsilon}
\varepsilon_N(t)  :=\displaystyle\sup_{u \le t} \left| \displaystyle\sum_{\mathbf{l}} \frac{\mathbf{l}}{N} \widetilde{Y}_{\mathbf{l}}\left(N\displaystyle \int_0^u \beta_{\mathbf{l}} \left(X^N(s)\right)ds\right) \right|.
\end{equation}
\vspace{0.1cm}

Therefore, using \eqref{F2}, \eqref{epsilon} and the Gronwall inequality, we can rewrite \eqref{distance1} as follows,

\begin{equation}
\label{distance2}
\begin{array}{lll}
|X^N(t) - X(t) | &\le& |X^N(0) - x_0 | + \varepsilon_N(t)  +  \displaystyle\int_0^t M|X^N(s) -X(s) |ds \\[0.4cm]
&\le& (|X^N(0) - x_0 | + \varepsilon_N(t) )e^{Mt},
\end{array}
\end{equation}
\noindent for all $t\ge0$. Thus, since $\varepsilon_N(t)$ is a non-decreasing function,

\begin{equation}
\label{distance3}
\sup_{s \le t} |X^N(s)- X(s)| \le  (|X^N(0) - x_0 | + \varepsilon_N(t) )e^{Mt}.
\end{equation}
\vspace{0.1cm}

To obtain \eqref{result1} we need the following lemma.

\medskip

\begin{lemma}
\label{lem1}
For every $t \ge 0$,

\begin{equation}
\displaystyle\lim_{N \to \infty} \varepsilon_N(t) = 0 \ \ a.s.
\end{equation}
\end{lemma}

\medskip

\begin{proof}
First, note that we have the following equivalence in distribution
\begin{equation}
\frac{\widetilde{Y}_{\mathbf{l}}(Nu)}{N}  \stackrel{D}{=} \displaystyle\sum_{i=1}^N \frac{Y^i_{\mathbf{l}}(u)}{N} - u,
\end{equation}
where $Y^i_{\mathbf{l}} (u)$ are independent Poisson processes with rate $u$. Then, the strong law of large number implies that

\begin{equation}
\label{poisson}
\displaystyle\lim_{N \to \infty} \sup_{u \le v} | N^{-1} \widetilde{Y}_{\mathbf{l}} (Nu) | = 0 \ \ a.s., \ \ v \ge 0.
\end{equation}
\vspace{0.1cm}

Furthermore, by definition of $\bar{\beta}_{\mathbf{l}} $ above, it follows that

\begin{equation}
\label{epsilon2}
\varepsilon_N(t)  \le \displaystyle\sum_{\mathbf{l}}  |\mathbf{l} | \displaystyle\sup_{u \le t} \left| \frac{\widetilde{Y}_{\mathbf{l}}\left(N u \bar{\beta}_{\mathbf{l}} \right)}{N}\right| .
\end{equation}
\vspace{0.1cm}

If $u\le t$, a basic property of Poisson processes states that,

\begin{equation}
Y_{\mathbf{l}}\left(N u \bar{\beta}_{\mathbf{l}} \right) \le Y_{\mathbf{l}}\left(N t \bar{\beta}_{\mathbf{l}} \right),
\end{equation}
\vspace{0.1cm}
for each $\mathbf{l} \in \mathcal{J}$. Then,

\begin{equation}
\label{termbyterm}
\begin{array}{ll}
|\mathbf{l} | \displaystyle\sup_{u \le t} \left| \frac{\widetilde{Y}_{\mathbf{l}}\left(N u \bar{\beta}_{\mathbf{l}} \right)}{N}\right| &\le \displaystyle\frac{ |\mathbf{l} |}{N} \left( Y_{\mathbf{l}}\left(N t \bar{\beta}_{\mathbf{l}} \right) + Nt \bar{\beta}_{\mathbf{l}}\right)\\[0.3cm]
&=  |\mathbf{l} |\displaystyle\frac{ Y_{\mathbf{l}}\left(N t \bar{\beta}_{\mathbf{l}} \right)}{N} + |\mathbf{l} | t \bar{\beta}_{\mathbf{l}},
\end{array}
\end{equation}
\vspace{0.2cm}
for all $N\ge1$ and each $\mathbf{l} \in \mathcal{J}$. Now, note that,

\begin{equation}
\label{poisson2}
 \displaystyle\sum_{\mathbf{l}}  |\mathbf{l} |\displaystyle\frac{ Y_{\mathbf{l}}\left(N t \bar{\beta}_{\mathbf{l}} \right)}{N} \stackrel{D}{=} \displaystyle\sum_{\mathbf{l}}  |\mathbf{l} |\displaystyle\sum_{i=1}^N \frac{Y^i_{\mathbf{l}}(t \bar{\beta}_{\mathbf{l}})}{N}= \displaystyle\sum_{i=1}^N\frac{1}{N}\left(\displaystyle\sum_{\mathbf{l}}  |\mathbf{l} | Y^i_{\mathbf{l}}(t \bar{\beta}_{\mathbf{l}})\right),
\end{equation}
\vspace{0.1cm}

\noindent where $Y^i_{\mathbf{l}} (u)$ are again independent Poisson processes with instensity $u$. Finally, by the law of large numbers applied to the independent random variables in bracket at r.h.s of \eqref{poisson2}, and by condition \eqref{betafinite},

\begin{equation}
\label{dominate}
\begin{array}{ll}
\displaystyle\lim_{N \to \infty} \displaystyle\sum_{\mathbf{l}}  |\mathbf{l} |\displaystyle\frac{ Y_{\mathbf{l}}\left(N t \bar{\beta}_{\mathbf{l}} \right)}{N} &=  \displaystyle\sum_{\mathbf{l}}  | {\mathbf{l}} | t \bar{\beta}_{\mathbf{l}}  \\[0.4cm]
&=  \displaystyle\sum_{\mathbf{l}} |\mathbf{l} | \displaystyle\lim_{N \to \infty} \displaystyle\sum_{i=1}^N \frac{Y^i_{\mathbf{l}}(t \bar{\beta}_{\mathbf{l}})}{N}.
\end{array}
\end{equation}

That is, by \eqref{termbyterm} we can interchange the limit and summation for the expression at r.h.s in \eqref{epsilon2}. Therefore, using \eqref{poisson}, the Lemma follows.
\end{proof}

\medskip

Finally, from inequality \eqref{distance3}, by condition $\lim_{N \to \infty} X^N ( 0 ) = x_0$, and Lemma \ref{lem1}, then for every $t \ge 0$,

\begin{equation}
\displaystyle\lim_{N \to \infty} \sup_{s \le t} |X^N(s)- X(s)| = 0 \ \ a.s.
\end{equation}
\end{proof}


\bigskip

Next section is devoted to our results about modeling a specific biological feedback loop, that is, our non-symmetric clock module. Clearly, the application of the Theorem \ref{th1} allows us to analyse the qualitative behaviour of the concentrations of the molecules involved in the interactions. Therefore, we are able to include the role of stochasticity in the gene expression, but we also simplify the qualitative study through a bifurcation analysis for the associated dynamical system.

\bigskip


\section{The associated dynamical system: formulation and bifurcation analysis}
\label{sec:bifur}

In this section we include a qualitative study of the dynamical behaviours of the non-symmetric clock module. We use the ideas of previous sections to analyse a dynamical system related to stochastic evolution of the concentration of molecules involved in our feedback loop.

First, we briefly summarize the ideas in previous sections. Initially, we stated a mean-field TDSIM $(\sigma_t)_{t \ge 0}$ on $\Omega_{\Lambda}$ and defined its dynamics by the rates of transition \eqref{taxas.ising}. Thus, the concentrations of each component in our non-symmetric loop (types $ A,  B$ and $ C $) were characterized by a density-profile process, defined in \eqref{densityp}. In the Ising case, the only independent variables are the densities of activated types that we denoted $X^N_i(t)$, $i\in \{A,B,C\}$. The rate of a jump in density-profile process was defined in \eqref{taxa} and \eqref{beta}, based on the rate of corresponding transition in TDSIM. After that, in Section \ref{sec:perfil}, we also studied the thermodynamic limit of the particle system ($N \to \infty$), proving that the density-profile process converges to deterministic trajectories governed by non-linear differential equations (see Theorem \ref{th1}).

Finally, in this section we show that depending on the parameter values, the magnetization random path can either converges to a unique stable fixed point (if $-1<J<2$), converges to one of a pair of stable fixed points (for $J<-1$), or asymptotically evolves close to a deterministic orbit in $\mathbb{R}^3$ (when $J>2$).

Therefore, we follow Theorem \ref{th1} to study the dynamics of our non-symmetric clock module by analysing the associated dynamical system. However, we remark that the behaviours of the dynamical system are deterministic approximations of the stochastic evolutions of the jump processes in our approach.

The system of differential equations associated to the cycle-interaction module in Figure \ref{fig:nonsym}, and with rates given by \eqref{taxas.ising} will be expressed by

\begin{equation}
\label{dinamico.bifurcacoes}
\dot{x}_i = (1-x_i) e^{ 2[-\delta Jx_{a(i)} - (1-\delta) J x_{h(i)} +  \kappa_i]} - x_i e^{-2[-\delta Jx_{a(i)} - (1-\delta) J x_{h(i)} +  \kappa_i]},
\end{equation}
\vspace{0.1cm}

\noindent for $i \in \{A, B, C\}$, $J \in \mathbb{R}$, $\delta \in [0,1]$, and $\kappa_i \in \mathbb{R}$. As usual, we will consider $\kappa_i = J/2$, $i \in \{A,B,C\}$, that it, proportional to the molecules interaction parameter. The particular value is to reduce the number of parameters and to obtain suitable fixed points in dynamical analysis (steady state of concentrations $(1/2,1/2,1/2)$).

\begin{figure}[ht]
\centering
\includegraphics[scale=0.28]{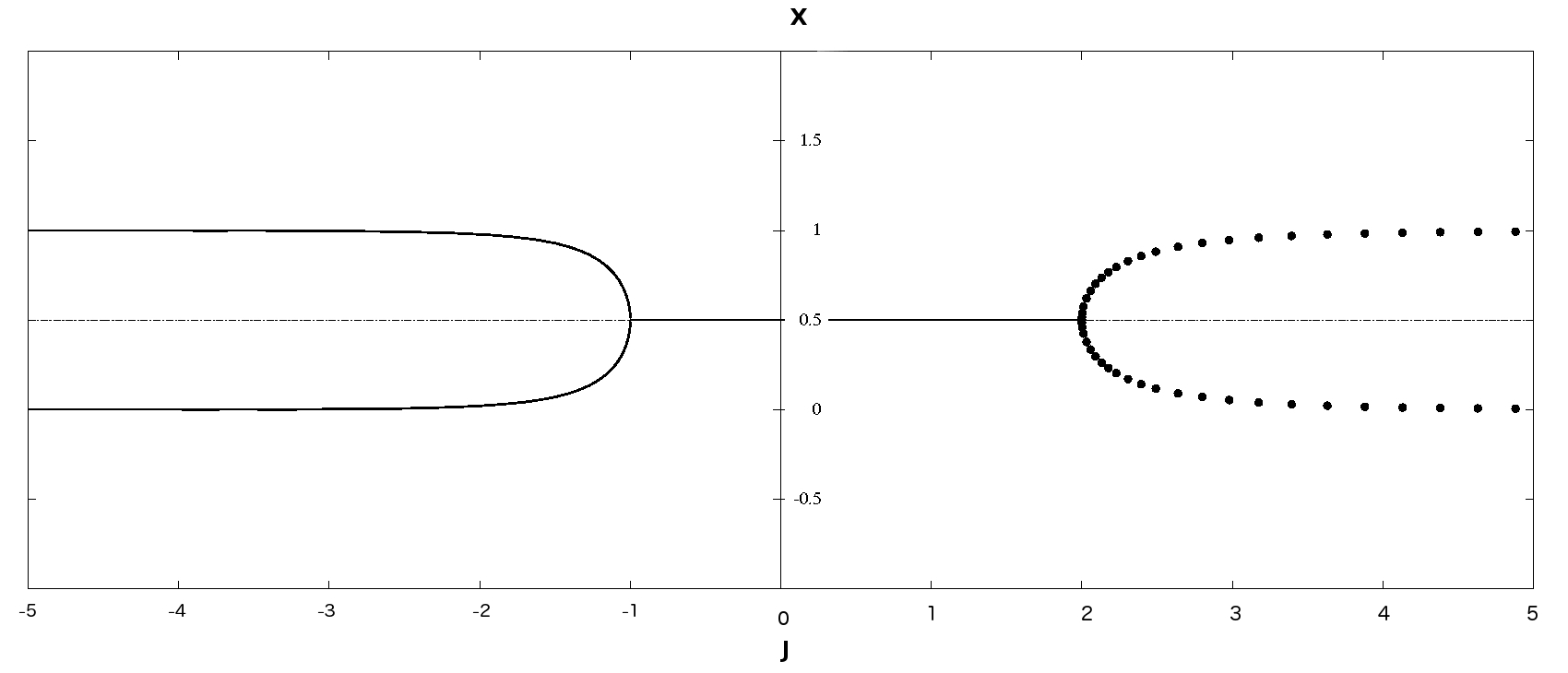}
\caption{{\scriptsize Bifurcation diagram for the non-symmetric clock module. The vertical axis represents the concentration of one of the three molecules. At right hand, we see a Hopf bifurcation with respect to the parameter $J$, that measures the strength of the interactions. Solid lines indicate stable points; dotted lines indicate unstable points; black circles indicate maximum and minimum values of stable orbits; bifurcation diagram obtained with xpp-aut \cite{Xppaut}.}}
\label{fig:bif}
\end{figure}

Our next result is the statement of a bifurcation analysis to guarantee for the non-symmetric clock module that: for $0 < J < J_c = 2$, the concentration of all three components ($A, B, C$) remains stable, but oscillates with large amplitude as soon as $J$ increases past the threshold $J_c=2$. For $J < 0$, we can see the appearance of two stable points when $J < -1$. These situations are showed in Figure \ref{fig:bif}.

\bigskip

\begin{proposition}
\label{bifurcacoes.relogio}
Consider the dynamical system \eqref{dinamico.bifurcacoes}, with $\delta \in D = [0,1/2) \cup (1/2,1]$, and $\kappa_i = J/2$, then
\ 
\begin{itemize}
\item [a)] For $J < 0$, there is a bifurcation at $J_c=-1$: the fixed point $(1/2, 1/2, 1/2)^T$ loses the stability and appear two stable points for $J < J_c$.
\item [b)] If $J > 0$, there is a Hopf (or Andronov-Hopf) bifurcation at $J_c=2$.
\end{itemize}
\end{proposition}

\medskip

\begin{proof}

First, it is easy to see that the dynamical system has a fixed point at $x^0 = (1/2, 1/2, 1/2)^T$, because $\dot{x} = F(x^0, J, \delta) = (0, 0, 0)^T$, for all $J \in \mathbb{R}$ and $\delta \in D$.

From now on, we consider the dynamical system $\dot{y} = F^*(y, J, \delta)$, where $y=(y_1, y_2, y_3)^T$ and $x_i = \frac{1}{2} +y_i$. Then, the fixed point for the new system is $\underline{0} = (0,0,0)^T$.

Following the theory in Perko \cite{Per} and Kuznetsov \cite{Kuz}, we could guarantee that near $\underline{0} = (0,0,0)^T$, the dynamical system is locally topologically equivalent to its linearisation $\dot{y} = A y$, where $A$ denotes the Jacobian matrix $\frac{dF^*}{dy}$ evaluated at $\underline{0} $. Then, $y^0$ is stable if all eigenvalues $\lambda$ of $A$ satisfy $Re \ \lambda < 0$, and stability is lost when one of the real parts becomes positive. 

 Notice that $A = \partial F^*(\underline{0})$ is given by the expression

\begin{equation}
\label{A2}
- 2 \left( \begin{array}{ccc}

1 &
 (1-\delta) J &
 \delta J
  \\
   \delta J&
1 &
(1-\delta) J
  \\
  (1-\delta) J &
   \delta J &
1
  \end{array} \right).
\end{equation}
\vspace{0.1cm}

The eigenvalues of this matrix are

\begin{equation}
\begin{array}{c}
\label{autovalores1}
\lambda_1 = - J - 1,\\[0.2cm]
\lambda_2 =  J(1  - 2\delta) \sqrt{3} i + (J-2), \text{ and}\\[0.2cm]
\lambda_3 = J (2  \delta - 1)\sqrt{3} i +(J- 2).
\end{array}
\end{equation}
\vspace{0.1cm}

Therefore, the fixed point $y^0$ loses the stability at $J_c= - 1$, when $\lambda_1$ crosses the imaginary axis through the origin. Thus, we have two stable fixed points. This bifurcation is called fold (or pitchfork) bifurcation. The values of these fixed points could be calculated by solving: $\sinh(2Jy) + 2y \cosh(2Jy) = 0$, thus we obtain the relation

\begin{equation}
\label{fixedpoint}
J = \frac{1}{4y} \cdot \log{\frac{1-2y}{1+2y}}
\end{equation}
\vspace{0.1cm}

\noindent  for $y\neq 0$. We notice the symmetry around $J$-axis, this fact could be seen at left hand of Figure \ref{fig:bif}.

On the other hand, the stability is lost when $\lambda_2$ and $\lambda_3$, symmetric around the real axis, cross the imaginary axis. This fact occurs when $J=2$. The phenomenon includes the appearance of a stable orbit, this kind of bifurcation is known as Hopf bifurcation.
\end{proof}

\bigskip


The oscillations experimentally verified in the repressilator (Elowitz and Leibler \cite{E&L}) are explained in this dynamical behaviour through the (supercritical) Hopf bifurcation with respect to the interaction strength parameter $J$.

\bigskip
\subsection{Deterministic macroscopic evolution in details}

Let us finally show a detailed analysis of the dynamical behaviours of the system using linearisation tools. In fact, we do not study the dynamical system \eqref{dinamico.bifurcacoes}, however we analyse the linear system that approximates it.

We will exhibit a qualitative analysis of the concentration of the molecules involved in the interactions. The idea is a better understanding of the influence of parameter $\delta$ in the dynamical behaviours.

Therefore, we must translate and rotate the coordinate system $X_1X_2X_3$ onto a new system $Z_1Z_2Z_3$. Particularly, this new coordinate system has the origin at $(1/2,1/2,1/2)$. Thus, $Z_3$ will be the diagonal, from $(0,0,0)$ to $(1,1,1)$; the $Z_2$-axis is directed to $x=(0,1,1/2)$; and the $Z_1$-axis is oriented to $x=(3/4,3/4,0)$.

Equivalently, remember that $y=(y_1, y_2, y_3)^T$ and $x_i = \frac{1}{2} +y_i$, for $i=\{1,2,3\}$. We can say that the system $Z_1Z_2Z_3$ is obtained by the following two operations in the system $Y_1Y_2Y_3$: 

\medskip

\begin{itemize}
\item [\textbf{1.-}] Rotate the $Y_1Y_2Y_3$-system around the $Y_3$-axis by $\alpha = 45^{\circ}$ (anticlockwise direction);

\item [\textbf{2.-}] Rotate the $Y_1Y_2Y_3$-system again about the new rotated $Y_2$-axis by $\beta = 55^{\circ}$ (clockwise direction);
%
\end{itemize}

\medskip

Finally, to obtain the new coordinates of the system, we have $(z_1, z_2, z_3) = (y_1, y_2, y_3) \cdot R$, with

\begin{equation}
\begin{array}{llc}
\label{R.rotacao}
R&= &  \begin{pmatrix}
 \frac{1}{\sqrt{6}}&  
- \frac{1}{\sqrt{2}}&
\frac{1}{\sqrt{3}}
  \\[0.2cm]
 \frac{1}{\sqrt{6}}&
\frac{1}{\sqrt{2}}&
\frac{1}{\sqrt{3}}
  \\[0.2cm]
- \frac{2}{\sqrt{6}}&
0&
\frac{1}{\sqrt{3}}
\end{pmatrix} .
\end{array}
\end{equation}
\vspace{0.2cm}

 Again, by the theory in \cite{Kuz,Per}, the behaviour of the dynamical system in $Z_1Z_2Z_3$, could be studied by linearisation, that is, we analyse

	\begin{equation}
\label{sistema.z}
\dot{z} = \textbf{A} z, \ \ \text{with} \ \ z = \begin{pmatrix} z_1\\ z_2 \\z_3 \end{pmatrix},
\end{equation}	
\vspace{0.2cm}

\noindent where $\textbf{A} = R^{-1} \cdot \partial F^*(\underline{0}) \cdot R$, because $R$ is an orthonormal matrix. Then,

\begin{equation}
\label{A.sistema.z}
\textbf{A} =  \begin{pmatrix}
J-2&  
\sqrt{3}J(2\delta -1)&
0
  \\[0.3cm]
-\sqrt{3}J(2\delta -1)&
J-2&
0
  \\[0.3cm]
0&
0&
-(2J+2)
\end{pmatrix} ,
\end{equation}
\vspace{0.2cm}

\noindent and expression \eqref{sistema.z} represents the following system:

\begin{equation}
\label{system.z}
\left\{\begin{array}{l}  \dot{z_1} = (J-2) z_1 + \sqrt{3}J(2\delta -1) z_2,\\[0.3cm]
 \dot{z_2} = - \sqrt{3}J(2\delta -1) z_1  + (J-2) z_2,\\[0.3cm]
 \dot{z_3} = -(2J+2) z_3. \end{array} \right.
\end{equation}
\vspace{0.2cm}

Notice that the behaviour over $Z_3$-axis is independent of the dynamics on $Z_1Z_2$-plane. Then, we could establish the following two elementary conclusions:

\medskip

\begin{itemize}
\item [\textbf{C1)}] For any $J > -1$, the trajectory on $Z_3$-axis goes to zero, because  $\dot{z_3}$ is negative for $z_3 > 0$, and positive for  $z_3 < 0$;

\item [\textbf{C2)}] For any $J < 2$, the dynamics on $Z_1Z_2$-plane goes to zero in both directions. In other words, it goes to the point (0,0). Of course, the velocity and direction depend on $\delta$ and $J$.
\end{itemize}

\medskip

Accordingly, given \textbf{C1)} and \textbf{C2)}, we conclude that for $-1 < J < 2$ the system has a stable point at $\underline{0} = (0, 0, 0)^T$.

Now, we study the dynamics for $J = - 1$. In this situation, $\dot{z_3} = 0$, for any $z_3$. Then, in adition with \textbf{C2)}, any initial condition $z(0)= (z_1(0), z_2(0), z_3(0))^T$ goes to $(0, 0, z_3(0))^T$, as $t \to \infty$. So, for $J < -1$ it is easy to see that, $\dot{z_3} < 0$ when $z_3$ is negative, and becomes positive when $z_3 > 0$. Namely, if the initial condition is established above the $Z_1Z_2$-plane, the trajectory goes to $(0,0,K)^T$, as $t \to \infty$, where $K\equiv K(J, \delta)$ is a positive constant that depends on $J$ and $\delta$. When the system starts at the bottom of the $Z_1Z_2$-plane, the trajectory goes to $(0,0,-K)^T$, as $t \to \infty$. These are the two stable points in Proposition \ref{bifurcacoes.relogio}.

Furthermore, considering $J = 2$, since \textbf{C1)}, we just need to study the dynamics over $Z_1Z_2$-plane. That is,

\begin{equation}
\label{system.z-2}
\left\{\begin{array}{l}  \dot{z_1} = 2\sqrt{3}(2\delta -1) z_2,\\[0.3cm]
 \dot{z_2} = - 2\sqrt{3}(2\delta -1) z_1.\end{array} \right.
\end{equation}
\vspace{0.2cm}

\begin{figure}[h!]
\centering
\includegraphics[scale=0.3]{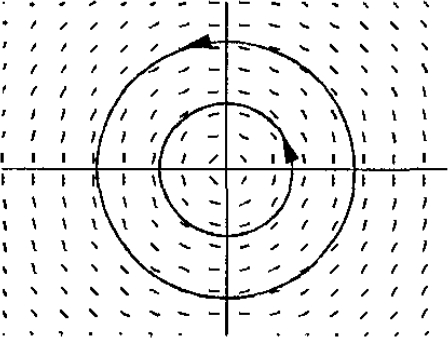}\hspace{0.7cm}
\includegraphics[scale=0.3]{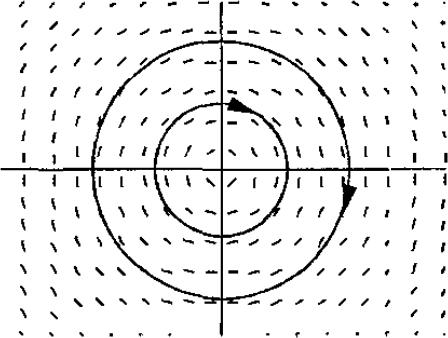}
\caption{{\scriptsize Phase portraits for $J=2$ and different values of $\delta$. At left hand, we see the behaviour for $\delta < 1/2$. The right picture show us the cycles for $\delta > 1/2$. The highest velocities will be reach at $\delta=0$ and $\delta=1$.}}
\label{fig:phaseportrait}
\end{figure}

We could change variables to polar coordinates. Let $z_1 = r cos\theta$, $z_2 = r sin\theta$. To derive a differential equation for $r$, we note $z_1^2+z_2^2 = r^2$, so $z_1\dot{z_1} + z_2\dot{z_2} = r\dot{r}$. Substituting for $\dot{z_1}$ and $\dot{z_2}$ yields $\dot{r} = 0$. Thus from $\dot{\theta} = (z_1\dot{z_2} -z_2\dot{z_1})/r^2$, we find $\dot{\theta} = - 2\sqrt{3}(2\delta -1)$. Therefore, the origin is a center, and any initial condition $z(0)= (z_1(0), z_2(0), z_3(0))^T$ evolves over the $Z_1Z_2$-plane on the circle with radius $r(0)= z_1(0)^2 + z_2(0)^2$, the direction and velocity of the cycle depend on $\delta$. The Figure \ref{fig:phaseportrait} shows us the behaviour of the cycle for different values of $\delta$.

\vspace{0.1cm}

\begin{figure}[h!]
\centering
\includegraphics[scale=0.3]{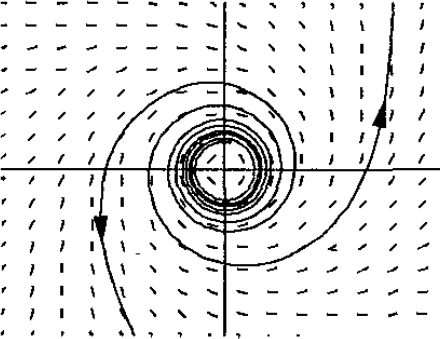}\hspace{0.7cm}
\includegraphics[scale=0.3]{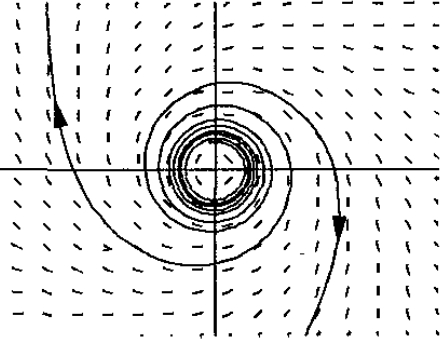}
\caption{{\scriptsize Phase portraits for $J >2$. The left picture shows us the behaviours for $\delta < 1/2$. At right side, we see the cycles for $\delta > 1/2$.}}
\label{fig:phaseportrait2}
\end{figure}

In the remaining case, for $J> 2$, we must study the two-dimensional linear system $\dot{z} = A z$ with

\begin{equation}
\label{phaseJ>4}
A=  \begin{pmatrix}
J-2&  
\sqrt{3}J(2\delta -1)  \\[0.3cm]
-\sqrt{3}J(2\delta -1)&
J-2
\end{pmatrix} 
\end{equation}
\vspace{0.1cm}

\noindent thus, changing variables to polar coordinates we obtain the system $\dot{r} = r (J-2)$, and $\dot{\theta} = -\sqrt{3}J(2\delta -1)$. Notice that, $\dot{r} >0$ for all $r>0$, and the value of $\dot{\theta}$ depends on $\delta$. See Figure \ref{fig:phaseportrait2}.

\bigskip


\section{Conclusions and comments} 
\label{sec:final}

On one hand, based on a particular feedback loop we have described a theoretical framework to model and analyse the non-linear dynamics of gene regulatory networks. In the literature there exist several different approaches to this aim. However, in previous sections we showed that our approach, based on the well studied Ising model, can include the inherent variability in gene expression, also this microscopic approach could capture detailed characteristics of the system interpreted as a promotion-inhibition circuitry. Thus, it may be particularly useful to illuminate how simple feedback loops manage to perform their basic functionalities and how these functionalities are further integrated into the whole cellular regulatory network.

On the other hand, the mathematical study of the type-dependent dynamics suggests a long-standing theme of statistical physics of nonequilibrium systems. That is, the question of the nature of their stationary measure, and in particular of the existence of phase transitions. Thus, it will be interesting for us to propose local interactions in type-dependent Ising models, instead of mean-field interaction, at this point we could consider several lattices (like $\mathbb{Z}^d$ or triangular lattice) and treat these questions.

Particularly, the study of the TDSIM with local interactions is closely related with works like Godr\`eche \cite{God1} and Godr\`eche and Bray \cite{God2} (see also de Oliveira \cite{MJO}). They considered a kind of asymmetric Ising dynamics, called directed Ising model. In \cite{God1,God2} was studied several dynamics, like Glauber or Metropolis. They obtained a large space of parameters defining the rates functions allowing irreversible Gibbsian Ising models, whenever the dynamics is not totally asymmetric.

Furthermore, bifurcations showed in the associated dynamical systems (Section \ref{sec:bifur}) suggest the existence of phase transitions in the type-dependent Ising model. Of course, it would also be natural to consider metastability issues \cite{OV} to characterize the dynamical behaviors of this particle system. In this sense, a related work is due to Koteck\'y and Olivieri \cite{K&O}, they studied a discrete-time Metropolis dynamics for a two dimensional ferromagnetic asymmetric Ising model, in which the vertical and horizontal interaction parameters have different values.

Finally, we say that another very interesting mathematical question is the characterization of loss and recovery of Gibbsianness in stochastic evolutions, considered in recent works, in the area of mathematical physics. For instance, van Enter et al. \cite{Fer} studied some Ising models for which this phenomenon occurs. In addition, Kulske and Le Ny \cite{K&L} initiated a fruitful research direction showing that Gibbs-non-Gibbs transitions can also be defined naturally for mean-field models.

\bigskip

\section*{Acknowledgments} The author thanks the anonymous referees for helpful suggestions. Prof. Eduardo J. Neves (IME/USP, Brazil) for many theoretical discussions. Prof. J. R. Mendon\c{c}a (EACH/USP, Brazil) and Alvaro Cerda (UFRO, Chile) for several comments. This work was successively supported by BecasChile, CONICYT and S\~ao Paulo Research Foundation FAPESP (Grant 2015/02801-6).

\bigskip


\end{document}